\documentclass[oneside, 10pt]{amsart}
\usepackage[cp1251]{inputenc}

\usepackage[T1]{fontenc}

\makeatletter

\providecommand{\LyX}{L\kern-.1667em\lower.25em\hbox{Y}\kern-.125emX\@}

 \theoremstyle{plain}
 \newtheorem{thm}{Theorem}[section]
 \numberwithin{equation}{section}
 \numberwithin{figure}{section} 
 \theoremstyle{plain}
 \theoremstyle{plain}
 \newtheorem{lem}[thm]{Lemma} 
 \theoremstyle{plain}
 \newtheorem{prop}[thm]{Proposition} 
 \theoremstyle{definition}
 \newtheorem{defn}[thm]{Definition}
 \theoremstyle{definition}
  
 \theoremstyle{remark}
 \newtheorem*{rem*}{Remark}

 \DeclareMathOperator*{\esssup}{ess\,sup}
\DeclareMathOperator{\supp}{supp}

\def\be{\begin{equation}\label}
\def\ee{\end{equation}}

\def\oH1{\buildrel\circ\over H\kern-.04in{}^1}

\def\supp{{\,supp\,}}

\def\bee{\begin{equation*}}
\def\eee{\end{equation*}}
\def\be{\begin{equation}}
\def\ee{\end{equation}}

\begin{document}

\title{Embedding Theorems and Boundary-value
 Problems for cusp domains.}

\author{V Gol'dshtein, M.Ju.Vasiltchik$^1$} \footnotetext[1]{The second author was partially supported by Russian Foundation for Basic Research (grant 06-01-00735)}

\address{$^{\textrm{*}}$Department of Mathematics, Ben Gurion University of
the Negev, Beer sheva, Israel \\
 vladimir@bgumail.bgu.ac.il
}

\address{$^{\textrm{**}}$Department of mathematics, Novosibirsk Technical University, Novosibirsk, Russia
}

\title{Embedding Theorems and Boundary-value
 Problems for cusp domains.}

\begin{abstract}
 We study the Robin boundary-value problem for bounded domains with isolated singularities. Because trace spaces of space $W^1_2 (D)$ on boundaries of such domains are weighted Sobolev spaces $L^{2, \xi}(\partial D)$  existence and uniqueness of corresponding Robin boundary-value problems depends on properties of embedding operators $I_1: W^{1}_2(D)\to L^{2}(D)$ and $I_{2}:W^{1}_2 (D)\rightarrow L^{2,\xi}(\partial D)$ i.e. on types of singularities.
We obtain an exact description of weights $\xi$ for bounded domains with 'outside peaks' on its boundaries. This result allows us  to formulate correctly  the corresponding Robin boundary-value problems for elliptic operators. Using compactness of embedding operators $I_1,I_2$, we prove also that these Robin boundary-value problems with the spectral parameter is of Fredholm type.

\end{abstract}

\maketitle

\section{Introduction}

An essentially self-contained presentation of a method for a correct
formulation and an investigation of the Robin boundary-value problem
for second order elliptic equations in domains with 'outside peaks'
 is given in this paper. This
study is initiated by works  \cite{r407}, \cite{r409}, \cite{GR},
\cite {GR1}. Quoted papers devoted to the Robin boundary-value
problem for Lipschitz domains and its strait forward generalizations.
It is well known that the embedding operator $I_1:W^{1}_2 (D)\to
L^{2}(D)$ is compact for Lipschitz domains \cite{M} and by \cite
{GR1} the operator $I_{2}:W^{1}_2(D)\rightarrow L^{2}(\partial D)$ is
also compact. Therefore the Robin boundary-value problem is of
Fredholm type for this class of domains \cite {GR1}.

For domains with H\"older type singularities the second embedding operator
$I_{2}:W^{1}_2(D)\rightarrow L^{2}(\partial D)$ does not necessarily exist, 
because traces of functions from $W^{1}_2(D)$ do not necessarily
belong to the space $L^{2}(\partial D)$. It means that
even a correct formulation of the Robin problem for such domains
depends on properties of a trace space of $W^{1}_2 (D)$. One of
possible descriptions of trace spaces for smooth bounded domains
with isolated singularities of an 'outside peak' type was proposed in  \cite{V1}, \cite{V2} by the second author.
For such domains the trace space for  $W^{1}_2 (D)$ does not necessarily
coincides with $L^2(D)$ and can be described with the
help of corresponding weights $\xi$ that depend on
singularity types. These results allow us to formulate correctly the Robin boundary-value problem  with the help of the weights $\xi$.

We will prove in this paper that embedding
operators $I_{2}:W^{1}_2(D)\rightarrow L^{2,\xi}(\partial D)$ are
compact for such weights. Because both embedding operators $I_1: W^{1}_2(D)\to L^{2}(D)$ and $I_{2}:W^{1}_2 (D)\rightarrow L^{2,\xi}(\partial D)$ are compact  the Robin boundary-value problem with the spectral parameter is of Fredholm type.

Elliptic boundary-value problems were studied in
numerous books and papers. We mention \cite{GT} and
\cite{LU}, where many references can be found. In \cite{M}
embedding theorems for a variety of  non-smooth domains have
been studied.

Main results of this paper were announced in \cite {VG}.

\section {Main results}

In this section we introduce necessary definitions and formulate main results. 

Recall the classical definition of Sobolev space $W^1_p(D)$.
Let $D$ be an open subset of $\mathbb R^n$, $n\geq 2$.  Define Sobolev space
$W^1_p(D)$, $1\leq p<\infty$ as a normed space of locally summable, weakly differentiable functions $f:D\to\mathbb R$ equipped with the following norm:
$$
\|f\mid W^1_p(D)\|=
\left[ \biggr(\int\limits_D|f|^p(x)\,dx\biggr)^{2/p}+ \biggr( \int\limits_D
|\nabla f|^p(x)\,dx\biggr)^{2/p} \right]^{1/2}.
$$
Here $\nabla f$ is the weak gradient of the function $f$.

\subsection{Weak Robin problem for smooth domains}

Let $G\in R^n$ be a bounded smooth domain i.e. its boundary $\partial G$ is a smooth compact $(n-1)$-dimensional manifold. Remember classical Robin boundary-value problem
\begin{equation}
\label{eq1}
\sum\limits_{i,j=1}^n \frac{\partial}{\partial x_i}\left(a_{ij}(x)\frac{\partial u}{\partial x_j}\right)+\sum\limits_{i=1}^n b_i(x) \frac{\partial u}{\partial x_i}+a(x)u=f(x),\ x\in G,
\end{equation}
\begin{equation}
\label{eq2}
\frac{\partial u}{\partial N}+\sigma(x)u=\mu(x),\ x\in\partial G.
\end{equation}
Here the functions $a_{i,j}(x)=a_{j,i},b_i(x),a(x),f(x)$
are smooth functions defined on $G$ and $\sigma(x),\mu(x)$ are
smooth functions defined on $\partial G$. 

As usually $\frac{\partial u}{\partial N}=\sum\limits_{i,j=1}^n a_{ij}
\frac{\partial u} {\partial x_j}\cos\left(\vec{n},x_i\right)$,
where $\vec{n}$ is the exterior unit normal vector at a point $x
\in \partial G$.

Suppose also that the following condition of the uniform ellipticity holds
\begin{equation}
\label{eq5}
 C_1|\xi|^2\leq \sum\limits_{i,j=1}^n
a_{ij}(x)\xi_i \xi_j \leq C_2|\xi|^2.
\end{equation}
for all $\xi=\left(\xi_1,\ldots,\xi_n\right)\in \mathbb{R}^n$.
Here $C_1$ and $C_2$ are positive constants.

This formulation is equivalent to the following weak formulation.
A function $u\in W^1_2 \left(G\right)$ is a weak
solution of the Robin problem if
\begin{equation}
\label{eq3}
\begin{split}
\sum\limits_{i,j=1}^n \int\limits_{G}a_{ij}(x) \frac{\partial u}{\partial x_j} \frac{\partial \overline{\eta}}{\partial x_i}dx -
\sum\limits_{i=1}^n \int\limits_{G} b_i(x) \frac{\partial u}{\partial x_i} \overline{\eta}dx-
\int\limits_{G}a u \overline{\eta} dx \\
=-\int\limits_G f \overline{\eta}dx-\int\limits_{\partial G}\sigma u \overline{\eta}dS_x+\int\limits_{\partial G}\mu\overline{\eta} dS_x.
\end{split}
\end{equation}
for all $\eta\in W^1_2 \left(G\right)$. Here $dS_x$ is the standard surface measure.

For the weak formulation smoothness of functions $f,\sigma,\mu$ is not necessary and can be replaced by the following weak conditions: $a_{i,j},b_i, a \in L^{\infty}(G)$, $f\in L^2\left(G\right)$, $\sigma\in L^{\infty}(\partial G)$,$\mu \in L^{2}(\partial G)$.

\subsection {Weak Robin problem for domains with an 'outside peak'}

Suppose domain $G$ is not smooth at one isolated point.

Assumptions on the functions $a_{ij}=a_{ji}$, $f$ and $a$ are the same for any such domain: $f\in L^2\left(G\right)$, $a \in L^{\infty}$ and $a_{ij}=a_{ji} \in L^{\infty}(G)$. Additional conditions for  $b_i$, \ $\sigma$,\ $\mu$ depends on properties of the boundary and will be formulated later only for domains with 'outside peaks' .

The next definition is a formal description of domains with an 'outside peak'.

\begin{defn} We call a bounded domain $G\subset \mathbb{R}^{n}$ a domain of class $OP_{\varphi}$ if

1. There exists such point $O\in\partial G$ that
 $\partial  G \setminus \{O\}$ is a smooth $(n-1)$-dimensional manifold of the class
$C^1$.

2. Let $\Omega \subset\mathbb{R}^{n-1}$ be a bounded domain
of the class $C^1$ and $\varphi\in C^1\left([0,1]\right)$ be a
smooth function such that $\varphi(0)=\varphi'(0)=0$ and $\varphi'(t)>0$ for
$t\in(0,1)$.Denote $x'=\left(x_1,\ldots,x_{n-1}\right)$.

There exists a neighborhood $U(O)$ of $O$ that can be represented as
\begin{equation}
\label{eq4}
  U(O)\cap G=\left\{x=(x',x_n)\in \mathbb{R}^n:\ 0<x_n<1,
\frac{x'}{\varphi(x_n)}\in\Omega\right\}.
\end{equation}
for an appropriate choice of a coordinate system with the origin $O$ in $\mathbb{R}^n$.
\end{defn}
 The point $O$ is a top of an 'outside peak'.  We will study problem
\eqref{eq3} for domains of class $OP_{\varphi}$.

Denote  $L^{p,\xi}\left(\partial G\right)$ such space of measurable
functions defined on $\partial G$ that
\begin{equation*}
 \int\limits_{\partial G} \left|f(x)\right|^p \xi(x) dS_x\equiv
\left\Vert f\right\Vert^p_{p,\xi,\partial G}< \infty.
\end{equation*}
Here $\xi: \partial G \rightarrow R$ is a fixed nonnegative
measurable function (a weight).

\emph{We are ready to specify all necessary assumptions  for Robin
boundary-value problem \eqref{eq3} in domains of
class $OP_{\varphi}$. }

Suppose domain $G$ belongs to the class $OP_{\varphi}$.

A function $u\in W^1_2 \left(G\right)$ is a weak
solution of the Robin problem if
\begin{equation}
\label{eq3a}
\begin{split}
\sum\limits_{i,j=1}^n \int\limits_{G}a_{ij}(x) \frac{\partial u}{\partial x_j} \frac{\partial \overline{\eta}}{\partial x_i}dx -
\sum\limits_{i=1}^n \int\limits_{G} b_i(x) \frac{\partial u}{\partial x_i} \overline{\eta}dx-
\int\limits_{G}a u \overline{\eta} dx \\
=-\int\limits_G f \overline{\eta}dx-\int\limits_{\partial G}\sigma u \overline{\eta}dS_x+\int\limits_{\partial G}\mu\overline{\eta} dS_x.
\end{split}
\end{equation}
for all $\eta\in W^1_2 \left(G\right)$. 

Here functions $a_{ij},b_i, a \in L^{\infty}(G)$, $f\in L^2\left(G\right)$.
The following uniform ellipticity condition holds
\begin{equation}
\label{eq5a}
 C_1|\xi|^2\leq \sum\limits_{i,j=1}^n
a_{ij}(x)\xi_i \xi_j \leq C_2|\xi|^2.
\end{equation}
for all $\xi=\left(\xi_1,\ldots,\xi_n\right)\in \mathbb{R}^n$ and some positive constants $C_1$ and $C_2$.

Additional assumptions for functions $ \sigma, \mu$ depends on function $\varphi$ and are essential in a neighborhood of singularity point $O \in \partial G$. Roughly speaking these assumptions must correlate with the exact description of the trace space of $W^1_2 (G)$ on the boundary $\partial G$. Reasons for the following assumptions will be made clear during proofs of the main results.

The function $\sigma:\partial G \rightarrow R$ satisfies to the following inequality:
\begin{equation}
\label{eq7}
\esssup_{x\in \partial G} \frac{|\sigma(x)|}{\varphi(x_n)}=M_\sigma<\infty.
\end{equation}
The function $\mu$ belongs to $L^{2,\frac{1}{\varphi}}\left(\partial G\right)$.
This is equivalent to $\frac{\mu}{\varphi} \in L^{2,\varphi}\left(\partial G\right)$.

Remember that functions $b_i$,\ $a$,\ $f$,\ $u$,\ $\eta$ (see
\eqref{eq3}) are complex valued functions.

\subsection{Main results}

Let $I_1$ be the embedding operator of $W^1_2 \left(G\right)$ into $L^2\left(G\right)$, and $I_2$ the embedding operator of $W^1_2 \left(G\right)$ into $L^{2,\varphi}\left(\partial G\right)$. By \cite{V1} the space $L^{2,\varphi}\left(\partial G\right)$ contains traces of $W^1_2 \left(G\right)$ on $\partial G$. Existence, boundedness and compactness of the operator $I_1$ is well known (see, for example \cite{M}). Existence and boundedness of the operator $I_2$ is proved in  \cite{V1}. Compactness of $I_2$ is the main technical goal of this paper.

Remember that $W^1_2 \left(G\right)$ with an inner product
\begin{equation}
\label{eq9}
\left[u,\eta\right]=\int\limits_G \left[\sum\limits_{i,j=1}^n  \frac{\partial u}{\partial x_i} \frac{\partial \overline{\eta}}{\partial x_j}+ u\overline{\eta}\right] dx
\end{equation}
is a Hilbert space.

We adopt a general well known statement of functional analysis to our study.
\begin{prop}
\label{Th1} Suppose $G \in OP_{\varphi}$  and operators $I_1$ and
$I_2$ are compact. Then the weak Robin problem \eqref{eq3a}-\eqref{eq7} is of Fredholm type, i.e. the problem can
be reduced to an operator equation on $W^1_2(G)$
\begin{equation}
\label{eq8}
\left(I+A\right)u=F,
\end{equation}
where $I$ is the identity and $A$ is a compact operator.
\end{prop}

Compactness of the embedding operator $I_2$ is a content of the following result:
\begin{thm}
\label{Th2}
If $G \in OP_{\varphi}$, then the embedding operator $I_2$ is compact.
\end{thm}

This theorem is a special case $p=2$ of the corresponding result for Sobolev spaces $W^1_p(G)$
that will be proved in Section 4.
Combining two previous results and using compactness of the embedding operator $I_1$ \cite{M} we obtain one of the main results of this study.
\begin{thm}
\label{Th3}
If $G \in OP_{\varphi}$, then the weak Robin problem \eqref{eq3a}-\eqref{eq7} is of Fredholm type.
\end{thm}
Proof of main results is based on the exact descriptions of  the trace spaces of Sobolev spaces $W_p^1\left(G\right)$  on boundaries of $OP_{\varphi}$-domains  \cite{V1}. For readers convenience we reproduce here this description.

Further the
relation $A\sim B$ means that a two sided
inequality $C_1<\frac{A}{B}<C_2$ with constants $0<C_1<C_2<\infty$
depending on $G$ only is correct.

Denote $TW_p^1\left(G\right)$
a normed space of $W_p^1\left(G\right)$-function traces on $\partial G$
with the quotient norm
\begin{equation}
\|f\|_{TW_p^1\left(G\right)}
=\inf\left\{\|F\|_{W_p^1\left(G\right)}: F\in W_p^1\left(G\right),
F\mid_{\partial G}= f\right\}.
\end{equation}
Let $E(x,y) = \max\{\varphi(x_n),\varphi(y_n)\}$ and $\sigma(x,y) =
\chi\left(\frac{|x - y|}{E(x,y)}\right)$, where $\chi$ is
 the indicator $[0,1]$.

 Let $G \in OP_{\varphi}$, $1<p<\infty$. By \cite{V1}
\begin{equation}
\left\Vert f\right\Vert_{TW_p^1} \sim \left\Vert f
\right\Vert_{p,\varphi,\partial G} +
\left\{\iint\limits_{G\times G}\frac{|f(x) -
f(y)|^p}{|x-y|^{n+p-2}}\sigma(x,y)dS_xdS_y\right\}^{1/p}.
\end{equation}
It means, in particular, that $f / \partial G \in L^{p,\varphi}$.

\section {Proof of Proposition \ref{Th1}}

The space $W^1_2 (G)$ is a Hilbert space with the inner product \eqref{eq9}. Because the condition \eqref{eq5} this inner product is equivalent to the following new one
\begin{equation}
\label{eq9a}
\left[u,\eta\right]=\int\limits_G \left[\sum\limits_{i,j=1}^n a_{ij}(x) \frac{\partial u}{\partial x_i} \frac{\partial \overline{\eta}}{\partial x_j}+ u\overline{\eta}\right] dx.
\end{equation}
Rewrite the weak Robin boundary problem in terms of this
inner product
\begin{equation}
\label{eq9b}
\left[u,\eta\right]-\sum\limits_{i=1}^n \int\limits_{G} b_i(x) \frac{\partial u}{\partial x_i} \overline{\eta}dx-
\int\limits_{G}(a+1) u \overline{\eta} dx +\int\limits_{\partial G}\sigma u \overline{\eta}dS_x\\
=-\int\limits_G f \overline{\eta}dx+\int\limits_{\partial G}\mu\overline{\eta} dS_x.
\end{equation}

Every integral in \eqref{eq9b} can be considered as
a complex valued linear functional on $W^1_2 \left(G\right)$:
\begin{equation*}
\begin{split}
l_{1,u}(\eta)&=-\sum\limits_{i=1}^n \int\limits_{G} b_i(x) \frac{\partial u}{\partial x_i}(x)\ \overline{\eta}(x)\ dx,\quad
l_{2,u}(\eta)=-\int\limits_{G} (a+1) u(x) \overline{\eta}(x) \ dx,\\
l_{3,u}(\eta)&=\int\limits_{\partial G} \sigma(x) u(x) \overline{\eta}(x) \ dS_x,\quad
l_4(\eta)=\int\limits_{\partial G} \mu(x) \overline{\eta}(x)\ dS_x,\\
l_5(\eta)&=-\int\limits_G f(x) \overline{\eta}(x) \ dx.
\end{split}
\end{equation*}

Because $a,b_i \in L^{\infty}(G)$,$f \in L^2(G)$, $u,\frac{\partial u}{\partial x_i}\in L^2(G)$ for any $i=1,2,...,n$ boundedness of functionals $l_{1,u}$,$l_{2,u}$,$l_5$ follows from Cauchy-Bunyakovski inequality.

Let us prove boundedness of the functional  $l_{3,u}$.
Using Cauchy-Bunyakovski inequality, condition \eqref{eq7}, the trace
description for domains of class $OP_{\varphi}$ \cite{V1} and boundedness of $I_2$ we obtain the
following inequality:
\begin{equation}
\label{eq11}
\begin{split}
\left|l_{3,u}(\eta)\right|= \left| \int\limits_{\partial G} \sigma
u \overline{\eta}\ dS_x \right| =\left| \int\limits_{\partial
G} \frac{\sigma}{\varphi} \varphi u \overline{\eta}\ dS_x
\right| \leq M_{\sigma} \left\Vert u\right\Vert_{2,\varphi,\partial
G}\cdot
\left\Vert \eta\right\Vert_{2,\varphi,\partial G} \\
\leq M_{\sigma}\left\|I_2\right\|^2 \left\Vert
u\right\Vert_{W^1_2 \left(G\right)}
\left\Vert\eta\right\Vert_{W^1_2 \left(G\right)}.
\end{split}
\end{equation}

 In \eqref{eq11} $W^1_2 \left(G\right)$-norm of functions is induced by the inner product\eqref{eq9}. Constants $M_{\sigma}$ and $\left\|I_2\right\|$ depends only on $G$ and function $\sigma$.
Therefore functional $l_{3,u}$ is bounded. A similar argument is correct for $l_4$:

\begin{equation}
\label{eq11a}
\begin{split}
\left|l_4(\eta)\right|= \left| \int\limits_{\partial G} \mu
\overline{\eta}\ dS_x \right|=\left| \int\limits_{\partial
G} \frac{\mu}{\varphi} \varphi  \overline{\eta}\ dS_x
\right| \leq \left\Vert \frac{\mu}{\varphi} \right\Vert_{L^{2,\varphi}(\partial
G)} \cdot
\left\Vert \eta\right\Vert_{H^1 (G)} \\
\leq \left\|I_2\right\| \left\Vert
\mu \right\Vert_{L^{2, \frac {1}{\varphi}} (\partial G) } \cdot
\left\Vert\eta\right\Vert_{H^1\left(G\right)}.
\end{split}
\end{equation}

We will use notations $C$,\ $C_1$,\ $C_2$,\ldots for different
positive constants. 

By the Riess
theorem there exist such bounded operators $B_i: W^1_2 \left(G\right)
\rightarrow W^1_2 \left(G\right)$,\ $i=1,2,3$ that
$
l_{i,u}(\eta):=\left[B_i u,\eta\right], \eta\in W^1_2 \left(G\right).
$

Denote $F:=B_4\mu+B_5 f$. Rewrite \eqref{eq3a} using  $B_i$ and $F$
\begin{equation*}
\left[u+(B_1+B_2+B_3)u,\eta\right]=\left[F,\eta\right],
\end{equation*}
We will use also a short notation  $A:=B_1+B_2+B_3$. 

Let us prove compactness of operator $A$. The operator
$A$ is compact if operators $B_1$,\ $B_2$,\ $B_3$ are compact.
It is enough to prove that for any operator $B_i$, $\,i=1,\,2,\,3$ an image of a weakly convergent sequence contains a strongly convergent subsequence.

Let a sequence $\{u_k\}_1^{\infty}$ weakly converges to $u_0$ in 
$W^1_2 (G)$. By continuity of $B_i$ the sequence
$B_iu_k$ weakly converges to $B_iu_0$. By compactness of operators $I_1$
and $I_2$ we can suppose (without loss of generality) that sequences $\{u_k\}_1^{\infty}$ and
$\{B_iu_k\}_1^{\infty}$ strongly converge in $L^2(G)$ and
$L^{2,\varphi}(\partial G)$ correspondingly. 

Let us start from $B_1$. For simplicity we will use short notations $l_i, i=1,2,3$ unstead of $l_{i,u_k-u_m}$. By the definition of $B_1$ 
\begin{gather}
   [B_1(u_k-u_m),\,B_1(u_k-u_m)]=\left| l_1(B_1(u_k-u_m))\notag\smallskip \right|\\
   \leq M \|B_1u_k-B_1u_m \|_{L^2(G)} .\notag
\end{gather}
Here $M$ is a positive constant.
Because $[v,\,v]\sim\|v\|^2_{W^1_2(G)}$ we have 
$$\|B_1u_k-B_1u_m\|_{W^1_2(G)}^2 \leq C\,\|B_1u_k-B_1u_m\|_{L^{2}(G)} $$
for a positive constant $C$.

The last inequality means that $\{B_1u_k\}_1^{\infty}$
strongly converges to $B_1u_0$  in $W^1_2 (G)$. Therefore operator $B_1$ is compact.

A similar argument is correct for operator $B_2$: 
\begin{gather}
   \|B_2(u_k-u_m)\|^2_{W^1_2(G)}\sim[B_2(u_k-u_m),\,B_2(u_k-u_m)]\notag\smallskip\\
   =\left| l_2(B_2(u_k-u_m)) \right| \leq C \,
   \|B_2u_k-B_2u_m\|_{L^2(G)} \notag
\end{gather}
for a positive constant $C$.
Therefore sequence $\{B_2u_k\}_1^{\infty}$ strongly converges to $B_2u_0$  in $W^1_2(G)$ and operator $B_2$ is compact.

For operator $B_3$ we will use compactness of $I_2$.
Using similar arguments we have
\begin{gather}
   \|B_3(u_k-u_m)\|^2_{W^1_2(G)}\sim
   [B_3(u_k-u_m),\,B_3(u_k-u_m)]\notag\smallskip\\
   =\left| l_3(B_3(u_k-u_m))\right| 
   \leq M_{\sigma}\,\|u_k-u_m\|_{L^{2,\varphi}(\partial
   G)}\cdot \|B_3(u_k-u_m)\|_{L^{2,\varphi} (\partial G)} \notag
\end{gather}
proving compactness  of $B_3$.

Hence operator $A=B_1+B_2+B_3$ is compact. Therefore the operator $I+A$ is of the Fredholm type. Proposition \ref{Th1} proved.

\section { Compactness of embedding operator $I_p:W^1_p(G) \rightarrow L^{p,\varphi}\left(\partial G\right)$}

Theorem \ref {Th2} is a special case for $p=2$ of the following result
\begin{thm}
\label{Th2p} If $G \in OP_{\varphi}$, then the embedding operator
$I_p:W^1_p(G) \rightarrow L^{p,\varphi}\left(\partial G\right)$ is
compact.
\end{thm}

Fix a domain $G \in OP_{\varphi}$. Without loss of generality we can suppose that
$0<\varphi'(x_n)<1$ for any $0\leq x_n \leq 1$.
can suppose that a top of peak $O$ coincides with the origin of
coordinates and there exists such a neighborhood $V(O)$ of $O$
that $V(O) \bigcap G$ can be represented as
\begin{equation}
\label{eq14}
 \left\{x=(x',x_n)\in \mathbb{R}^n: 0< x_n <1, \mid x'\mid <
 \varphi(x_n)\right\}.
 \end{equation}
Here notice that $x=(x',x_n)\in \mathbb{R}^n$,\ $x'\in\mathbb{R}^{n-1}$,\ $x_n\in\mathbb{R}$. Without 

For the proof of Theorem \ref{Th2p} we need the following technical result:  
\begin{lem}
\label{Lem}

Any function  $F\in W^1_p\left(G\right)$,\ $1<p<\infty$, has the following representation
\begin{equation*}
F(x)=\alpha(x_n)+R(x),
\end{equation*}
where $\alpha\in W^1_p\left(G\right)$ is a function of one
variable, $R\in W^1_p\left(G\right)$, and the following
inequality
is correct
\begin{equation}
\label{eq13}
\left\Vert\varphi^{-1}(x_n) R\right\Vert_{L^p\left(G\right)}\leq
C \left\Vert F \right\Vert_{W^1_p\left(G\right)}<\infty,
\end{equation}
for some constant $C$ that depends only on $G$.
\end{lem}
For readers convenience we give here a complete, comparatively simple and independent proof of this result that was proved by the second author in  \cite{V1}.

 {\bf Proof}

  Choose nonnegative functions $h\in
C^\infty_0\left(\mathbb{R}^{n-1}\right)$ and $K\in
C^\infty_0\left(\mathbb{R} \right)$ such that\\
$\supp h\subset
B_1(0'):=\left\{x'=\left(x_1,\ldots,x_{n-1}\right)\in
\mathbb{R}^{n-1}:\ |x'|<1 \right\}$,\\ $\supp K\subset
\left[\frac{1}{2},1\right]$,\ $\int\limits_{\mathbb{R}^{n-1}}
h(x') dx'=1$,\ $\int\limits_{-\infty}^{\infty} K(t) dt=1$. 
Suppose
\begin{equation*}
\Omega(x,r,\theta',y)=\varphi^{-n}(x_n)   h\left(\frac{x'+r\theta'}{\varphi(x_n)}\right)r^{n-2}K\left(\frac{y}{\varphi(x_n)}\right),
\end{equation*}
where $r\geq 0$,\ $y\in\mathbb{R}^1$,\ $\theta'\in
S_1(0')=\partial B_1(0')$. We will use a family of kernels

$$\Omega_{\varepsilon}:=\Omega\left(x,\frac{r}{\varepsilon},\theta',\frac{y}{\varepsilon}\right).$$

Denote
\begin{equation*}
F_{\varepsilon}(x)=\int\limits_{S_1(0')}dS_{\theta'} \int\limits_{0}^{\infty} dr
\int\limits_{-\infty}^{\infty} F(x'+r\theta',x_n+y) \Omega\left(x,\frac{r}{\varepsilon},\theta',\frac{y}{\varepsilon}\right) \frac{dy}{\varepsilon^2}.
\end{equation*}
The proof of this Lemma can be divided onto four parts. \\
1. Let us demonstrate that 
 $F_{\epsilon} \rightarrow F$ in
$L^{p,loc} \left(G\right)$ for $\epsilon \rightarrow 0$. Using the change of
variables $r\vartheta' = z', y = z_n$ we obtain that 
\begin{equation*}
F_\varepsilon (x) =\int\limits_{R^n} F(x+z)h\left(\frac{\varepsilon
x' + z'} {\varepsilon
\varphi(x_n)}\right)K\left(\frac{z_n}{\varepsilon\varphi(x_n)}\right)\frac{dz}{\varepsilon^n
\varphi^n (x_n)}.
\end{equation*}\

Let us prove that $F_\varepsilon (x)\rightarrow F(x)$ for
$\varepsilon\rightarrow 0$ in $L^{p,loc} ( G)$. Let $U$ be
a compact subset of $G$. Denote $x^0_n = \inf\{x_n : x =
(x',x_n) \in U\}$ Then $x^0_n > 0$. Using monotonicity of  $\varphi$
and compactness of $U$ we have that  $ \varphi(x^0_n)< \varphi(x_n) < C_0
\varphi(x^0_n)$, where $x = (x',x_n) \in U$ and the constant $C_0$ depends only on $U$ 
and $\varphi$. By construction of $K$ we have  $\frac 1 2
\leq \frac{z_n}{\varepsilon \varphi(x_n)} \leq 1$. Hence
\begin{equation*}
K\left(\frac{z_n}{\varepsilon \varphi(x_n)}\right) \leq C_1
\chi\left(\frac{z_n}{C_0 \varepsilon \varphi(x^0_n)}\right)
\end{equation*}
Because $h\left(\frac{\varepsilon x' + z'}{\varepsilon
\varphi(x_n)}\right)\neq 0$ for $\left|\frac{\varepsilon x'
+z'}{\varepsilon \varphi(x_n)}\right| < 1$, then
\begin{equation*}
\left| h\left(\frac{\varepsilon x' + z'}{\varepsilon
\varphi(x_n)}\right)\right| \leq C_2 \chi\left(\frac{|z'|}{2C_0
\varepsilon \varphi(x^0_n)}\right).
\end{equation*}
In these inequalities $C_1 = \sup\{K(\tau): \tau \in \mathbf{R}\}$,
$C_2 = \sup\{|h(u')|: u' \in \mathbf{R^{n-1}}\}$ and $\chi$ 
is the indicator of $[0;1]$. The second inequality is based on 
the inequality  $ |z'| \leq \varepsilon |x'| + \varepsilon
\varphi(x_n) \leq 2\varepsilon \varphi(x_n) \leq 2C_0 \varepsilon
\varphi(x^0_n)$ that is a direct consequence of an obvious inequality
$|\varepsilon x' + z'| \leq \varepsilon
\varphi(x_n)$. 

Because
\begin{equation*}
\int_{\mathbf{R^n}} h\left(\frac{\varepsilon x' + z'}{\varepsilon
\varphi(x_n)}\right) K\left(\frac{z_n}{\varepsilon
\varphi(x_n)}\right) \frac{dz}{\varepsilon^n \varphi^n (x_n)} = 1,
\end{equation*}
we get
\begin{equation*}
|F_\varepsilon (x) - F(x) | \leq C\int_{\mathbf{R^n}} |F(x+z) - F(x)
| \chi\left(\frac{z_n}{C_0 \varepsilon \varphi(x^0_n)}\right)
\chi\left(\frac{|z'|}{2C_0 \varepsilon \varphi(x^0_n)}\right)
\frac{dz}{\varepsilon^n \varphi^n (x^0_n)}.
\end{equation*}
Integrating both parts of the previous inequality in degree $p$ on  $U$ 
and using Minkovsky inequality 
\begin{multline*}
\|F_\varepsilon - F\|_{L^p (U)} \leq \\ C  \int_{R^n}  \left( \int_U \mid F(x+z)
-F(x)\mid^p dx \right)^\frac{1}{p} \chi\left(\frac{z_n}{C_0 \varepsilon
\varphi(x^0_n)}\right) \chi\left(\frac{\mid z'\mid}{2C_0 \varepsilon
\varphi(x^0_n)}\right) \frac{dz}{\varepsilon^n \varphi^n (x^0_n)}
\\ \leq C \sup\{\|F(x+z) - F(x)\|_{L^p (U)}: |z| \leq 2C_0 \varepsilon
\varphi(x^0_n)\}.
\end{multline*}
Using continuity of function $F \in L^p(G)$ in the sense of $L^p$ (see \cite{S}, ch.1)
 the proof can be finished. 

Therefore we obtained the following integral representation 
\begin{equation*}
F(x)=F_1(x)-\int\limits_{0}^{1}\frac{\partial}{\partial\varepsilon} F_{\varepsilon}(x) d\varepsilon.
\end{equation*}
for any $F \in L^p(G)$.
 
 We will use short notation
\begin{equation*}
\quad R(x)=-\int\limits_{0}^{1}\frac{\partial}{\partial\varepsilon} F_{\varepsilon}(x) d\varepsilon.
\end{equation*}

 Recall that
    \begin{equation*}
    F_1 (x)= \int\limits_{-\infty}^\infty dy \int\limits_S
    dS_\theta' \int\limits_0^\infty F(x'+r\theta',x_n
    +y)\Omega(x,r,\theta',y) dr.
    \end{equation*}
   
 Using the following change of variables $x'+r\theta'=z', x_n +y=z_n, r^{n-2}
   dr d\theta'=dz', dy=dz_n$ we obtain

   \begin{equation}
   \label{aeq1}
   \alpha(x_n): = F_1(x)=
   \int\limits_{R^n}F(z)h\left( \frac {z'}{\varphi(x_n)}\right) K\left( \frac{z_n
   -x_n}{\varphi(x_n)} \right) \frac{dz}{\varphi^n(x_n)}.\
   \end{equation}

   Here we used a new notation $\alpha(x_n): = F_1(x)$ emphasizing that function
   $F_1(x)$ depends only on $x_n$.

   2. Let us demonstrate that $\alpha\in W_p^1 (G)$ for $1< p <\infty$.
   By construction of $\alpha$ we have
    \begin{multline}
    \label{aeq2}
    \parallel\alpha\parallel_{L^p(G)}^p = \int\limits_G
    \mid\alpha(x_n)\mid^p dx = \int\limits_0^1 dx_n
    \int\limits_{\mid x'\mid<\varphi(x_n)} \mid\alpha(x_n)\mid^p dx' \\ \leq
    \int\limits_0^1 \int\limits_{\mid x'\mid<\varphi(x_n)}
    \left|\int\limits_{R^n}
    \mid F(z)\mid\cdot  h(\frac{z'}{\varphi(x_n)}) \cdot  K(\frac{z_n
    -x_n}{\varphi(x_n)})  \frac{dz}{\varphi^n(x_n)} \right|^p dx'dx_n\
    \end{multline}

    Remember that $\supp h \subset B_1(0')$ and $ \supp K
    \subset[1/2;1]$. Therefore the domain of integration for $z$
    is the following set $U:=\{z: \mid z'\mid<\varphi(x_n), x_n + \varphi(x_n)/2 < z_n <
    x_n + \varphi(x_n)\}$ which belongs to $G$.

     There exist $\delta > 0$  and $o<c<x_n$ such that
    \begin{equation*}
    \varphi(x_n +\delta\varphi(x_n)) = \varphi(x_n) +
    \varphi'(c)\delta\varphi(x_n) = \varphi(x_n)(1 + \delta\varphi'(c))
    \end{equation*}
    By monotonicity of $\varphi$ we have for $\delta= 1/2$ the following inequality
    \begin{equation}
    \label{aeq3}
    \frac{\varphi(x_n)}{2}\leq\varphi(z_n)\leq\frac{3\varphi(x_n)}{2}
    \end{equation}

    The integrand in \eqref{aeq2} can be rewritten in more convenient way
    \begin{equation*}
    \varphi^{-n} \mid F\mid\cdot\mid h\mid\cdot\mid K\mid =
    \mid F \mid\cdot(\mid h\mid^{\frac 1{p}}\mid K\mid^{\frac 1{p}}\varphi^{-\frac n{p}})\cdot(\mid h\mid^{\frac{p-1}{p}}\mid K\mid^{\frac{p-1}{p}}\varphi^{-\frac{(p-1)n}{p}}).
    \end{equation*}
    Using H\"older inequality and changing order of integration we obtain
    \begin{equation}
    \label{F} \left\Vert \alpha\right\Vert _{L_{p}(G)}=C\:\int_{G}\left|F(z)\right|^{p}\Psi(z,x_{n})dz
    \end{equation}
    where
    \[
C=\left\{ \int_{G}\Psi(z,x_{n})dz\right\} ^{p-1},\]
\[
\Psi(z,x_{n})=\int_{0}^{1}dx_{n}\int_{\left\{ \left|x'\right|<\varphi\left(x_{n}\right)\right\} }h\left(\frac{z'}{\varphi\left(x_{n}\right)}\right)K\left(\frac{z_{n}-x_{n}}{\varphi\left(x_{n}\right)}\right)\frac{dx'}{\varphi^{n}\left(x_{n}\right)}.\]
We will estimate $\Psi(z,x_{n})$. 
    
    Let us change first $x_n$ onto
$\tau:=\frac{z_n-x_n}{\varphi(x_n)}$.
Hence
\begin{equation*}
\frac{dx_n}{\varphi(x_n)}=-\frac{d\tau}
{1+\frac{z_n-x_n(\tau)}{\varphi(x_n(\tau))}\ \varphi'(x_n)}.
\end{equation*}
Therefore
\begin{equation*}
\int\limits_{0}^{1} K\left(\frac{z_n-x_n}{\varphi(x_n)}\right)
\frac{dx_n}{\varphi(x_n)}\leq C \int\limits_{-\infty}^{\infty}
K(\tau) \frac{d\tau}{1-\varphi'(x_n(\tau))}< \infty,
\end{equation*}
because of the boundedness of function $K$ with compact support. Using the
same argument for \eqref{F} we obtain finally the following estimate
for the $L^p$-norm of $\alpha$
    \begin{equation*}
    \parallel\alpha\parallel_{L^p (G)}^p \leq C\int\limits_G
    \mid F(z)\mid^p dz.
    \end{equation*}
    Let us prove that $\alpha\in W^1_p(G)$.
    
Using the following change of variable in \eqref{aeq1}

    \begin{equation*}
    z' = \varphi(x_n)s',  z_n = x_n+ \varphi(x_n)s_n.
    \end{equation*}
   we get
    \begin{equation*}
    \alpha(x_n) = \int\limits_{R^n} F(\varphi(x_n)s',x_n +
    \varphi(x_n)s_n) h(s')K(s_n)ds.
    \end{equation*}
After differentiation we obtain

    \begin{multline*}
    \alpha'(x_n) =
    \int\limits_{R^n}(\varphi'(\nabla'F(\varphi(x_n)s',x_n+\varphi(x_n)s_n),s')+ \\
    D_n F(\varphi(x_n)s',x_n+ \varphi(x_n)s_n)(1+\varphi'(x_n)s_n)) h(s')K(s_n)
    ds.
    \end{multline*}
   Returning to the original variable we get finally
    \begin{multline*}
    \mid\alpha'(x_n)\mid \leq C  \int\limits_{R^n}
    \mid\nabla'F(z)\mid\cdot  h(\frac{z'}{\varphi(x_n)}) \cdot  K(\frac{z_n
    -x_n}{\varphi(x_n)}) \frac{dz}{\varphi^n (x_n)} \\+
    C \int\limits_{R^n}
    \mid D_nF(z)\mid\cdot  h(\frac{z'}{\varphi(x_n)}) \cdot  K(\frac{z_n
    -x_n}{\varphi(x_n)}) \frac{dz}{\varphi^n (x_n)}.
    \end{multline*}
    This inequality permits us to estimate $\parallel\alpha'\parallel_{L^p (G)}$
    by the same way as
    $\parallel\alpha\parallel_{L^p (G)}.$

    Therefore $\alpha \in W^1_p(G)$.

 3. Function $R$ can be represented as $R=R_1+R_2$ where
\begin{equation*}
R_1(x)=-\int\limits_0^1 \frac{d\varepsilon}{\varepsilon^{n+1}} \int\limits_{\mathbb{R}^n}
\left(\nabla' F(z), z'-x'\right)\ h\left(\frac{z'-(1-\varepsilon)x'}{\varepsilon\varphi(x_n)}\right)
K\left(\frac{z_n-x_n}{\varepsilon\varphi(x_n)}\right) \frac{dz}{\varphi^n(x_n)},
\end{equation*}
\begin{equation*}
R_2(x)=-\int\limits_0^1 \frac{d\varepsilon}{\varepsilon^{n+1}} \int\limits_{\mathbb{R}^n}
D_n F(z)\cdot (z_n-x_n)\
K\left(\frac{z_n-x_n}{\varepsilon\varphi(x_n)}\right)
 h\left(\frac{z'-(1-\varepsilon)x'}{\varepsilon\varphi(x_n)}\right)
 \frac{dz}{\varphi^n(x_n)}.
\end{equation*}
Here
 $\nabla'=\left(\frac{\partial}{\partial z_1},\ldots,\frac{\partial}{\partial z_{n-1}}\right)$,\ $z=(z',z_n)$,\ $D_n=\frac{\partial}{\partial z_n}$.

Let us verify this representation.
Remember first that
\[
    \Omega \left ( x,\frac r \epsilon,\theta',\frac y \epsilon \right ) =
    \frac{r^{n-2}}{\varphi^n (x_n)\epsilon^{n-2}}
    h \left ( \frac{\epsilon x'+r\theta'}{\epsilon\varphi(x_n)}\right ) K \left ( \frac y
    {\epsilon\varphi(x_n)} \right ).
\]
    Hence
    \begin{equation}
    \label{aeq4}
\frac{\partial}{\partial\varepsilon}\left[\varepsilon^{-2}\Omega\left(x,\frac{r}{\varepsilon},\theta',\frac{y}{\varepsilon}\right)\right]=
\end{equation}
\[
-\frac{nr^{n-2}}{\varepsilon^{n+1}\varphi^{n}\left(x_{n}\right)}hK-\frac{r^{n-1}}{\varepsilon^{n+2}\varphi^{n+1}\left(x_{n}\right)}\sum_{j=1}^{n-1}D_{j}h\theta_{j}-\frac{yr^{n-2}}{\varepsilon^{n+2}\varphi^{n+1}\left(x_{n}\right)}hK'=\]
\[
-\frac{K}{\varepsilon^{n+1}\varphi^{n}\left(x_{n}\right)}\:\frac{\partial}{\partial r}\left[r^{n-1}h\right]-\frac{r^{n-2}}{\varepsilon^{n+1}\varphi^{n}\left(x_{n}\right)}h\frac{\partial}{\partial y}\left[yK\left(\frac{y}{\varepsilon\varphi\left(x_{n}\right)}\right)\right]\]

    Using \ref{aeq4} we get the following expression for $R(x)$
   \[
    R(x) = -\int\limits_0^1
    \frac{\partial F_\epsilon}{\partial\epsilon} d\epsilon = 
    \]
    \[
    -\int\limits_0^1 d\epsilon\int\limits_{-\infty}^{\infty}
    dy\int\limits_S dS_\theta' \int\limits_0^{\infty}
    F(x'+r\theta',x_n
    +y)\frac{\partial}{\partial\epsilon}\left[\Omega(x',\frac
    r{\epsilon},\theta',\frac y{\epsilon})\epsilon^{-2}\right] dr=
    \]
    \[
-\int_{0}^{1}d\varepsilon\int_{S}dS_{\theta'}\int_{0}^{\infty}\left(\sum_{j=1}^{n-1}D_{j}F\theta_{j}\right)r^{n-1}h\left(\frac{\varepsilon x'+r\theta'}{\varepsilon\varphi\left(x_{n}\right)}\right)dr-\]
\[
\int_{0}^{1}d\varepsilon\int_{S}dS_{\theta'}\int_{0}^{\infty}h\left(\frac{\varepsilon x'+r\theta'}{\varepsilon\varphi\left(x_{n}\right)}\right)\frac{r^{n-2}}{\varepsilon^{n+1}\varphi^{n}\left(x_{n}\right)}dr\int_{-\infty}^{\infty}y\, D_{n}F\, K\left(\frac{y}{\varepsilon\varphi\left(x_{n}\right)}\right)dy. \]

   Using the following change of variables
    \[
    z' = x' + r\theta', dz' = r^{n-2} drdS_\theta', x_n+y=z_n,
    dz_n = dy.
    \]
     we obtain the claimed before decomposition $R=R_1+R_2$:
    \[
    R(x) =
    -\int\limits_0^1\frac{d\epsilon}{\epsilon^{n+1}\varphi^n
    (x_n)}\int\limits_{R^n} (\nabla'F(z),z'-x')K\left(\frac{z_n
-x_n}{\epsilon\varphi(x_n)}\right)h\left(\frac{z'-(1-\epsilon)x'}{\epsilon\varphi(x_n)}\right)
   dz -
   \]
   \[
    \int\limits_0^1\frac{d\epsilon}{\epsilon^{n+1}\varphi^n
    (x_n)}\int\limits_{R^n} D_n F(z)\cdot(z_n -x_n)K\left(\frac{z_n
-x_n}{\epsilon\varphi(x_n)}\right)h\left(\frac{z'-(1-\epsilon)x'}{\epsilon\varphi(x_n)}\right)dz
    \]
    \[
    = R_1 + R_2.
   \]

 4.  Similarly to the plane case  \cite{V2} and to the proof of
\eqref{eq13} we will demonstrate that $\varphi^{-1}R_i\in
L^p\left(G\right)$,\ $i=1,2$. We will prove this fact only for
$R_1$. For $R_2$ the proof is similar.

We have
\begin{multline}
\label{eq15}
 \left|\varphi^{-1}(x_n)R_1(x)\right| \leq  \\
\int\limits_0^1 d\varepsilon \int\limits_{\mathbb{R}^n} |\nabla'
F(z)|\cdot\left|\frac{z'-x'}{\varepsilon\varphi(x_n)}\right|\
h\left(\frac{z'-(1-\varepsilon)x'}{\varepsilon\varphi(x_n)}\right)
K\left(\frac{z_n-x_n}{\varepsilon\varphi(x_n)}\right)
\frac{dz}{\varepsilon^{n}\varphi^n(x_n)}.
\end{multline}

 Denote
by $X(x)$ the extension of $\left|\varphi^{-1}(x_n)R_1(x)\right|$
on $\mathbb{R}^n\setminus G$ by zero and by $Y(z)$ the
extension of $|\nabla' F(z)|$ by zero on $\mathbb{R}^n\setminus
G$. 
Recall generalized Minkovski inequality. Let $\psi(x,y)$be a measurable
nonegative function defined on $A\times B$ where $A\in R^{n}$, $B\in R^{m}$
are measurable sets. Then for any $1\leq p<\infty$\[
\left\{ \int_{A}\left(\int_{B}\psi(x,y)dy\right)^{p}dx\right\} ^{1/p}\leq\int_{B}\left(\int_{A}\psi^{p}(x,y)dx\right)^{1/p}dy.\]

Using generalized Minkovski inequlity we get
\begin{equation}
\label{eq16}
\left\Vert X\right\Vert _{L_{p}(\mathbf{R}^{n})}\leq
\end{equation}
\[
\int_{0}^{1}\left\{ \int_{\mathbf{R}^{n}}\left[\frac{1}{\varepsilon^{n}\varphi^{n}\left(x_{n}\right)}\int_{\mathbf{R}}\Phi(z,x)\left(\int_{\mathbf{R}^{n-1}}h\left(\frac{z'-(1-\varepsilon)x'}{\varepsilon\varphi\left(x_{n}\right)}\right)dz'\right)^{\frac{p-1}{p}}Kdz_{n}\right]^{p}dx\right\} ^{1/p}d\varepsilon.\]
Here\[
\Phi(z,x)=\left\{ \int_{\mathbf{R}^{n-1}}Y^{p}(z)h\left(\frac{z'-(1-\varepsilon)x'}{\varepsilon\varphi\left(x_{n}\right)}\right)dz'\right\} ^{1/p}.\]
Because

\[
\left(\int_{\mathbf{R}^{n-1}}h\left(\frac{z'-(1-\varepsilon)x'}{\varepsilon\varphi\left(x_{n}\right)}\right)dz'\right)^{\frac{p-1}{p}}\leq C\left[\varepsilon\varphi\left(x_{n}\right)\right]^{\frac{(n-1)(p-1)}{p}}\]
we obtain finally
\begin{equation}
\label{new}
\left\Vert X\right\Vert _{L_{p}(\mathbf{R}^{n})}\leq C\int_{0}^{1}\left\{ \int_{\mathbf{R}^{n}}\left[\int_{\mathbf{R}}\Phi(z,x)K\left(\frac{z_{n}-x_{n}'}{\varepsilon\varphi\left(x_{n}\right)}\right)dz_{n}\right]^{p}\frac{dx}{\left(\varepsilon\varphi\left(x_{n}\right)\right)^{n+p-1}}\right\} ^{1/p}d\varepsilon.
\end{equation}

Using for the term in square brackets from the right hand side of
\eqref{new} H\"older inequality for variable  $z_n$ and taking
into account the inequality
\begin{equation*}
\left\{ \int\limits_{-\infty}^{\infty}
K\left(\frac{z_n-x_n}{\varepsilon\varphi(x_n)}\right) dz_n
\right\}^{p-1}\leq C
\left(\varepsilon\varphi(x_n)\right)^{p-1},
\end{equation*}
we obtain the following $L^p$-estimate
\begin{multline}
\label{eq17}
 \left\Vert X\right\Vert_{L^p\left(G\right)}\leq \\
C\int\limits_0^1 \left\{ \int\limits_{\mathbb{R}^n} Y^p(z) \left[
\int\limits_{\mathbb{R}^n}
h\left(\frac{z'-(1-\varepsilon)x'}{\varepsilon\varphi(x_n)}\right)
K\left(\frac{z_n-x_n}{\varepsilon\varphi(x_n)}\right)
\frac{dx}{(\varepsilon\varphi(x_n))^{n}} \right] dz
\right\}^{\frac{1}{p}} d\varepsilon.
\end{multline}
Denote by $J$ the term in square brackets on the right hand side of
\eqref{eq17}. By direct calculations we have
\begin{multline*}
J=\int\limits_{\mathbb{R}^n}
h\left(\frac{z'-(1-\varepsilon)x'}{\varepsilon\varphi(x_n)}\right)
K\left(\frac{z_n-x_n}{\varepsilon\varphi(x_n)}\right)
\frac{dx}{(\varepsilon\varphi(x_n))^{n}}= \\
\int\limits_{-\infty}^{\infty} K\left(\frac{z_n-x_n}{\varepsilon\varphi(x_n)}\right)
\left( \int\limits_{\mathbb{R}^{n-1}}
h\left(\frac{z'-(1-\varepsilon)x'}{\varepsilon\varphi(x_n)}\right) dx'
\right)
\frac{dx_n}{(\varepsilon\varphi(x_n))^{n}}\leq \\  C
\int\limits_{-\infty}^{\infty} K\left(\frac{z_n-x_n}{\varepsilon\varphi(x_n)}\right)
\frac{dx_n}{\varepsilon\varphi(x_n)}.
\end{multline*}

For the last inequality we used the following estimate\[
\int_{\mathbf{R}^{n-1}}h\left(\frac{z'-(1-\varepsilon)x'}{\varepsilon\varphi\left(x_{n}\right)}\right)dx'\leq C\left[\varepsilon\varphi\left(x_{n}\right)\right]^{n-1}.\]
This estimate can be proved by the following way. Let \[
A=\left\{ x'\in\mathbf{R}^{n-1}:\left|\frac{z'-x'}{\varepsilon\varphi\left(x_{n}\right)}\right|<1\right\} ,\; B=\left\{ x'\in\mathbf{R}^{n-1}:\left|z'-x'\right|<2\varepsilon\varphi\left(x_{n}\right)\right\} .\]

Using definition of class $OP_{\varphi}$ and inequalities\[
\left|\frac{z'-x'}{\varepsilon\varphi\left(x_{n}\right)}\right|-1\leq\left|\frac{z'-x'}{\varepsilon\varphi\left(x_{n}\right)}\right|-\frac{\left|x'\right|}{\varphi\left(x_{n}\right)}\leq\left|\frac{z'-(1-\varepsilon)x'}{\varepsilon\varphi\left(x_{n}\right)}\right|<1\]
follows that $A\subset B$. Therefore\[
\int_{\mathbf{R}^{n-1}}h\left(\frac{z'-(1-\varepsilon)x'}{\varepsilon\varphi\left(x_{n}\right)}\right)dx'\leq\int_{A}hdx'\leq\int_{B}hdx'\]

Let us change variable $x_n$ to
$\tau:=\frac{z_n-x_n}{\varepsilon\varphi(x_n)}=\tau$\
($x_n\rightarrow\tau$). Hence
\begin{equation*}
\frac{dx}{\varepsilon\varphi(x_n)}=-\frac{d\tau}
{1+\frac{z_n-x_n(\tau)}{\varepsilon\varphi(x_n(\tau))}\ \varphi'(x_n)}.
\end{equation*}
Therefore
\begin{equation*}
\int\limits_{-\infty}^{\infty} K\left(\frac{z_n-x_n}{\varepsilon\varphi(x_n)}\right)
\frac{dx}{\varepsilon\varphi(x_n)}\leq C \int\limits_{-\infty}^{\infty}
K(\tau) \frac{d\tau}{1-\varepsilon\varphi'(x_n(\tau))}\leq C,
\end{equation*}
where constant $C$ depends only on $F(x)$. Using the last
inequality \eqref{eq17} we obtain finally
\begin{equation}
\label{eq18}
 \left\Vert X\right\Vert_{L^p\left(\mathbb{R}^n\right)}\leq C\left\Vert Y\right\Vert_{L^p\left(\mathbb{R}^n\right)}\leq C\left\Vert F\right\Vert_{W^1_p\left(G\right)}.
\end{equation}
The Lemma is proved.

Using the previous Lemma \ref{Lem} we will prove the theorem \ref{Th2p}

{\bf Proof of Theorem \ref{Th2p}}
Let $\mathfrak{M}\subset W_p^1\left(G\right)$ be a bounded set and $\Lambda=I_p(\mathfrak{M})\subset L^{p,\varphi}\left(\partial G\right)$. To prove compactness of $I_2: W^1_p(G)\rightarrow L^{p,\varphi}\left(\partial G\right)$ we need to construct for any  $\varepsilon>0$ a finite $\varepsilon$-network of $\Lambda$.

By Lemma \ref{Lem}
any function $F\in\mathfrak{M}$ can be represented as $F=\alpha+R$,\ $\alpha=\alpha(x_n)$. Denote $\Lambda_{\alpha}$ the set of
all functions $\alpha=F-R$ such
that $F\in\mathfrak{M}$ and denote $\Lambda_R$ the set of all functions $R$ that
correspond to $F\in\mathfrak{M}$. Hence $\Lambda=\Lambda_{\alpha}+\Lambda_R$.

Because the embedding operator $I_1:W^1_p \rightarrow L^p(G)$ is compact the set
$\Lambda_{\alpha}$ is relatively compact in $L^p\left(G\right)$. Therefore
for any $\varepsilon_1>0$ there exists an $\varepsilon_1$-network
$\beta_1,\ldots,\beta_N$ of $\Lambda_{\alpha}$ in
$L^p\left(G\right)$. Denote $G_{\delta}=B_{\delta}(0)\cap G$ where
$B_{\delta}(0)=\left\{x\in\mathbb{R}^n:\ |x|<\delta\right\}$. Further we suppose that $0<\delta<1$. Using continuity of function $F \in L^p(G)$ in the sense of $L^p$ (see \cite{S}, ch.1),
we can found $\delta>0$ such that for any $j=1,2,\ldots,N$ the
following inequality is correct
\begin{equation*}
\left\Vert \beta_j\right\Vert_{p,G_{\delta}}<\varepsilon_1.
\end{equation*}
Hence for any $\alpha\in\Lambda_{\alpha}$ there exists such $\beta_j$ that
\begin{equation}
\label{eq19} \left\Vert \alpha\right\Vert_{p,G_{\delta}}\leq
\left\Vert \alpha-\beta_j\right\Vert_{p,G_{\delta}}+\left\Vert
\beta_j\right\Vert_{p,G_{\delta}}<2\varepsilon_1.
\end{equation}
By construction of the domain $G$ in a neighborhood of its peak
with the top $O$ (see \eqref{eq4}) we have by direct calculation
\begin{equation*}
\int\limits_{G_{\delta}} |\alpha(x_n)|^p
dx=\omega_{n-1}\int\limits_{0}^{\delta} |\alpha(x_n)|^p
\varphi^{n-1}(x_n) dx_n,
\end{equation*}
where $\omega_{n-1}$ is volume of the unit $n-1$-dimensional ball.

On the other hand we will demonstrate (using \eqref{eq4}) that

\begin{equation*}
\int\limits_{\partial G_{\delta}} |\alpha(x_n)|^p \varphi(x_n)
dS_x\sim \int\limits_{0}^{\delta} |\alpha(x_n)|^p
\varphi^{n-1}(x_n) dx_n.
\end{equation*}

 The argument is simple. The manifold $\partial G_{\delta}\setminus \{0\}$ is a finite union of  charts of the following type
\begin{multline}
\label{eq20}
\partial G_{{\delta},{\sigma}}^{(n-1)}=\biggl\{x=\left(x'',x_{n-1},x_n\right):\ 0<x_n<\delta,\
x_{n-1} \\ =\sqrt{\varphi^2(x_n)-|x''|^2},
|x''|=\left(\sum\limits_{i=1}^{n-2} x_i^2\right)^{1/2}<\sigma\varphi(x_n)\biggr\},
\end{multline}
where $\sigma\in (0,1)$. By direct calculations
\begin{equation*}
dS_x=\varphi(x_n)\sqrt{\frac{1+\left(\varphi'(x_n)\right)^2}{\varphi^2(x_n)-|x''|^2}} dx'' dx_n,
\end{equation*}
and by  the inequality $|x''|<\sigma\varphi(x_n)$
\begin{equation*}
\varphi(x_n)\sqrt{\frac{1+\left(\varphi'(x_n)\right)^2}{\varphi^2(x_n)-|x''|^2}}\sim 1.
\end{equation*}
Therefore
\begin{multline}
\int\limits_{\partial G_{\delta}^{(n-1)}} |\alpha(x_n)|^p
\varphi(x_n) dS_x\sim \int\limits_{0}^{\delta} dx_n
\int\limits_{\left\{|x''|<\sigma\varphi(x_n)\right\}}
|\alpha(x_n)|^p \varphi(x_n) dx'' \\ \sim \int\limits_{0}^{\delta}
|\alpha(x_n)|^p \varphi^{n-1}(x_n) dx_n.
\end{multline}
Combining this estimate and \eqref{eq19} we can conclude that integrals
\begin{equation*}
\int\limits_{\partial G_{\delta}} |\alpha(x_n)|^p  \varphi(x_n)
dS_x,\ \alpha\in\Lambda_{\alpha},
\end{equation*}
are uniformly small for all
$\alpha\in\Lambda_\alpha$ if $\delta >0$ is small enough.

To estimate the function $R \in \Lambda_R$ that corresponds to $F \in \mathfrak{M}$ we will use Lemma \ref{Lem}.  Denote $\partial_0 G = \{x\in\partial G:
    \mid x'\mid = \varphi(x_n), 0 < x_n < 1\}$. First we will estimate the integral
\begin{equation*}
\int\limits_{\partial_{0}G_{\delta}} \mid R(x)\mid^p \varphi(x_n)dS_x,
\end{equation*}
 for some $\delta > 0$.

 The boundary $\partial_{0} G_\delta$ can be covered by a finite number of charts
$\partial G_{\delta,\sigma}^{(k)}$,$k=1,2,...,n-1$, (see (\ref{eq20})). Let us prove estimates for $R$ for the chart $\partial G_{\delta,\sigma}^{(n-1)}$.
We will use notation $x=:(x",x_{n-1},x_n)\in\partial_0 G_{\delta,\sigma}^{(n-1)}$.Then
$$x_{n-1}=\sqrt{\varphi^2(x_n)-|x"|^2}, |x"|<\sigma\varphi(x_n).$$
Choose $s: 0<s<\sqrt{\varphi^2 (x_n) - \mid x"\mid^2}$. Then
\begin{equation*}
R(x) =R(x",s,x_n) + \int\limits_s^{x_{n-1}} D_{n-1}
R(x",\tau,x_n)d\tau
\end{equation*}
where $D_{n-1} = \frac{\partial}{\partial x_{n-1}}$.

Integrating in $s$ we obtain
\begin{multline}
\label{eq21}
\mid R(x)\mid^p \varphi(x_n) \leq C_1
\int\limits_0^{\sqrt{\varphi^2 (x_n) - \mid x"\mid^2}}
\mid R(x",s,x_n)\mid^p ds +  \\ 
C_2 \varphi(x_n)\left(\int\limits_0^{\sqrt{\varphi^2 (x_n) - \mid x"\mid^2}} \mid D_{n-1}
R(x",\tau,x_n)\mid d\tau\right)^p
\end{multline}
For the estimate \eqref{eq21} we used that $\sqrt{\varphi^2 (x_n)-|x"|^2}\sim
\varphi(x_n)$ if $\mid x"\mid < \sigma\varphi(x_n)$. Integrating
\eqref{eq21} on $\partial_0 G_{\delta,\sigma}^{(n-1)}$, taking into account
 $dS_x \sim dx"dx_n$ and using H\"older inequality for variable $\tau$ we get
the following inequality
\begin{equation*}
\begin{split}
&\int\limits_{\partial_0 G_{\delta,\sigma}^{(n-1)}}
\mid R(x)\mid^p\varphi(x_n)dS_x \leq  \\& C_1 \int\limits_0^\delta dx_n
\int\limits_{{\mid x"\mid<\sigma\varphi(x_n)}}
dx"\int\limits_0^{\sqrt{\varphi^2 (x_n) - \mid x"\mid^2}}
\mid R(x",s,x_n)\mid^p ds\\&+ C_2 \int\limits_0^\delta\varphi^p
(x_n)dx_n \int\limits_{{\mid x"\mid<\sigma\varphi(x_n)}}
dx"\int\limits_0^{\sqrt{\varphi^2 (x_n) - \mid x"\mid^2}} \mid D_{n-1}
R(x",\tau,x_n)\mid^p d\tau\\&\leq
C_1\int\limits_{G_\delta} \mid R(x)\mid^p dx +
C_2\varphi^p (\delta)\int\limits_{G_\delta }
\mid D_{n-1} R(x)\mid^pdx.
\end{split}
\end{equation*}
Combining with \eqref{eq18} we obtain
\begin{multline}
\int\limits_{\partial_0 G_{\delta,\sigma}^{(n-1)}}
\mid R(x)\mid^p\varphi(x_n)dS_x \leq  \\ C\varphi^p
(\delta)\left[\int\limits_G \left | \frac{R(x)}{\varphi(x_n)}\right |^p dx +
\int\limits_G \mid D_{n-1} R(x)\mid^p dx\right]
\end{multline}
Using this estimate and Lemma \ref{Lem} we get finally
\begin{equation}
\label{eq22}
 \left\Vert R \right\Vert_{p,\varphi,\partial_0 G_\delta}
\leq C\varphi(\delta)\left\Vert F\right\Vert_{W_p^1 (G)}.
\end{equation}

This estimate and corresponding estimates for integrals of
$\mid\alpha(x_n)\mid^p$ on $\partial_0 G_\delta$ for  $\alpha\in
\Lambda_\alpha$ allow us to conclude that for any $\epsilon_2 > 0$ there exists such $\delta_0 > 0$ that for all $\delta\in(0,\delta
)$ and for all $F\in\mathfrak{M}$ the following inequality is correct
\begin{equation}
\label{eq23}
 \left\Vert F\right\Vert_{p,\varphi,\partial_0 G_{\delta}} <
\epsilon_2.
\end{equation}
Let $\delta$ be an arbitrary number from (0,1). Denote
 $k_\delta$ a function of the class $C^1 ([0,1])$ with the following properties $0\leq k_\delta \leq 1$ on
$[0,1]$, $k_\delta(t) = 1$ when $\delta < t \leq 1$, $k_\delta (t) =
0$ when $0 \leq t \leq \delta/2$. Denote by
$\mathfrak{M}(\delta)$ traces of functions $F$ from $\mathfrak{M}$
on $G\backslash\overline{G_{\delta/2}}$. The set
$\mathfrak{M}(\delta)$ is bounded in $W_p^1
(G\backslash\overline{G_{\delta/2}})$
and therefore \cite{S} is compactly embedded into $L^{p,\varphi} (\partial_0
G\backslash\partial_0 G_{\delta/2}) \sim L^p (\partial_0 G
\backslash \partial_0 G_{\delta/2} )$ (for fixed
$\delta$). Let $\epsilon > 0$ be an arbitrary number. Choose a finite $\epsilon/2$-network for $\Lambda(\delta) =
I_2(\mathfrak{M}(\delta) ) \subset L^p(\partial_0 G \backslash
\partial_0 G_{\delta/2} )$ that is a finite set of functions $\tilde
\nu_1,\ldots,\tilde\nu_r$. Then  the set of functions $\nu_j = k_\delta
\tilde\nu_j, j=1,2,\ldots,r$ represents a finite $\epsilon/2$-network
for the set $k_\delta\Lambda(\delta) = \left\{ k_\delta F: F\in
\mathfrak{M}(\delta)\right\}$ into $L^{p,\varphi} (\partial_0 G)$.
It is possible to suppose that for the same $\delta$ the inequality
 \eqref{eq23} is correct with $\epsilon_2 = \epsilon/2$.

 Finally we will prove that
 the set of functions $\nu_1,\ldots,\nu_r$ represents an $\epsilon$-network of $\Lambda$ into $L^{p,\varphi} (\partial_0 G)$. The argument is standard. Let $F\in \mathfrak{M}$. Then
$F = (1 - k_\delta )F + k_\delta F$. Because $k_\delta F \in
k_\delta \Lambda(\delta)$ there exists such $\nu_i$ that
$\left\Vert k_\delta F - \nu_i\right\Vert_{p,\varphi,\partial_0 G} <
\epsilon/2$. Then by \eqref{eq23})
\begin{equation*}
\left\Vert F - \nu_i\right\Vert_{p,\varphi,\partial_0 G} \leq
\left\Vert(1 - k_\delta )F\right\Vert_{p,\varphi,\partial_0 G} +
\left\Vert k_\delta F - \nu_i\right\Vert_{p,\varphi,\partial_0 G} <
\epsilon/2 + \epsilon/2 = \epsilon.
\end{equation*}

Therefore $\nu_1,\ldots,\nu_r$ is a desire
$\epsilon$-network for $\Lambda = I_2 (\mathfrak{M})$ in $L^{p,\varphi}
(\partial_0 G)$.

Theorem \ref{Th2p} proved.


\begin{thebibliography}{50}

\bibitem{A}Sh.Agmon, Lectures on Elliptic Boundary Value Problems. D.Van
Nostrand Co. Princeton,
Toronto, New York, London, 1965.

\bibitem{BIN} O.V.Becov, V.P.Il'in, S.M.Hikolskii, Integral representation
of functions and embedding theorems (russian), Nauka, 1975.

\bibitem{CG} R.Courant, D.Hilbert, Methoden der mathematischen
Physik, vol 2, Springer Verlag, Berlin, 1937.

\bibitem{GT} D.Gilbarg, N. Trudinger, Elliptic partial differential
equations of second  order, Springer, Berlin, 2001
\bibitem{E} W.Evans, D.Harris, Sobolev embeddings for
generalized ridged domains, Proc. Lond. Math. Soc., 54, N3,
(1987), 141-175.

\bibitem{GhV} F.W.Gehring and J.Vaisala, Hausdorff dimension and
quasiconformal mappings.J.
London Math.Soc. (2),6, (1973),504-521.

\bibitem{GV}~V.M. Gol'dshtein, S.K. Vodop'janov, Prolongement des functions de classe L(2,1)
et applications quasiconformes.C.R.Ac.Sc., Paris, 1980, 290, 10, A453-456.

\bibitem{GR}V.M.Gol'dshtein, Yu.G. Reshetnyak, Quasiconformal Mappings and
Sobolev Spaces,
Kluwer Academic Publishers. Dordrecht, Boston, London. 1990.

\bibitem{GGR}V. Gol'dshtein, L.Gurov, A.Romanov, Homeomorphisms that induce monomorphisms
of Sobolev spaces, Isr. Journ. of Math., 91,1995,31-60.

\bibitem{GG}V.Gol'dshtein and L.Gurov. Applications of change of variable operators for
exact embedding theorems, Integr. Equat. Oper. Th., 19, (1994), 1-24.

\bibitem{GR1} V.Goldshtein, A.G.Ramm, Embedding operators for rough
domains, Math Ineq. and Applications, 4, N1, (2001) 127-141.

\bibitem{GR2} V.Goldshtein, A.G.Ramm, Embedding operators and boundary-value problems for
rough domains, IJAMM, 1, (2005) 51-72.

\bibitem{LU}  O.A. Ladyzhenskaya, N.N. Ural'tseva,
Linear and quasilinear elliptic equations,
Academic Press, New York, 1968.

\bibitem{M}V.Maz'ya, Sobolev Spaces, Springer Verlag, Berlin, 1985.

\bibitem{r407} A.G.Ramm, Inverse problems, Springer, New York, 2004.

\bibitem{r409} A.G.Ramm, M.Sammartino, Existence and uniqueness of
the scattering solutions in the exterior of rough domains,
in the book "Operator theory and its applications", Amer. Math. Soc.,
Fields Institute Communications, vol.25,  pp.457-472,

\bibitem{S} S.L.Sobolev, Some applicationns of functional analysis to
mathematical physics. Leningrad, Leningrad State Unuversity, 1950.

\bibitem{V1} M.Ju.Vasiltchik, Traces of functions from Sobolev space $W^1_p$ for domains with non Lipschitz boundaries. In: Modern problems of geometry and analysis. Novosibirsk, Nauka, (1989), 9-45.

\bibitem{V2} M.Ju.Vasiltchik, Necessary and sufficient conditions on traces of functions from Sobolev spaces for a plane domain with non Lipschitz boundaries. In: Studies on mathematical analysis and Riemannian geometry. Novosibirsk, Nauka, (1992), 5-29.

\bibitem{VG} M.Ju.Vasiltchik, V.M.Gol'dshtein, About solvability of third boundary-value problem for domains with peak. (Russian) Matematicheskie zametki, 78, 3, (2005), 466-468.

\bibitem{Z}W.P.Ziemer, Weakly Differentiable Functions, Springer Verlag, 1989




 \bibitem{Gris} Grisvard~P.
Elliptic problems in nonsmooth domains. ---~Pitman; Boston, 1985.
---~P.~410.
\end{thebibliography}
\end{document}